\documentclass[10pt,reqno]{article}
\usepackage{amsmath, amsthm, amssymb, stmaryrd}
\usepackage{hyperref}
\usepackage{enumerate}
\usepackage{url}
\usepackage{color}
\usepackage{tikz}
\usepackage{xcolor}
\usepackage{cancel}
\usetikzlibrary{patterns}
\usepackage[T1,T2A]{fontenc}
\usepackage[utf8]{inputenc}
\usepackage[all]{xy}
\usepackage{upgreek}
\usepackage{amsmath}

\theoremstyle{definition}

\theoremstyle{remark}

\def \le {\leqslant}
\def \ge {\geqslant}

\usepackage{tikz} 
\usetikzlibrary{calc,arrows,decorations.pathreplacing,fadings,3d,positioning}

\begin{document}

\begin{Large}
\centerline{On Furstenberg's Diophantine result
}
\end{Large}
\vskip+0.5cm

\begin{large}
\centerline{Dmitry Gayfulin and Nikolay Moshchevitin}
\end{large}

\vskip+1cm

{\bf 1. Introduction.}

\vskip+0.3cm

Let $a, b\ge 2$ be  coprime integer  numbers and
$$\Sigma = \{ a^ub^v:\,\, u,v\in \mathbb{Z}_+\}.
 $$
In 1967 Furstenberg \cite{F} using theory of disjointness of topological dynamical systems proved that for any irrational number 
$\alpha\in \mathbb{R}$ the sequence of fractional parts
$$
{\Sigma_\alpha=\biggl\{ \{q\alpha \} :\,\,\, q \in \Sigma \biggr\}}
$$
is dense in $[0,1]$. A short and elementary proof of this result was given by Boshernitzan \cite{B}. Both proofs were ineffective. An effective version of Furstenberg's Diophantine result was given by Bourgain, Lindenstrauss, Michel and Venkatesh in \cite{BLMV}. 
In particular, they proved the following two theorems.

\vskip+0.3cm
{\bf Theorem A.}\,
{\it
Let $a, b$ be multiplicatively independent integers i.e. $\frac{\log a}{\log b}\in\mathbb{R}\setminus\mathbb{Q}$. Suppose 
$\alpha \in \mathbb{R}$
is irrational and Diophantine-generic, that is there  exists $ k_1,k_2>0 $  such that 
$\left|\alpha -\frac{p}{q}\right| \ge k_1q^{-k_2}$ for all rationals $\frac{p}{q}$. 
Then
the set 
$$
{\biggl\{ \{a^ub^v\alpha \} :\,\,\, u,v \le M\biggr\}}
$$
is $\frac{1}{{(\log \log M)}^{\kappa(a,b)}}$-dense in $[0,1]$ for some effective positive $\kappa(a,b)$  
 and for $M\ge  M_0=M_0(k_1,k_2,a,b)$ where $M_0$ is also an effective constant.
}

\vskip+0.3cm

For our purpose the sets 
$${
\Sigma  (M)= \{ q\in \Sigma:\,\, q\le M\}}
\,\,\,\,\,
\text{and}
\,\,\,\,\,{\Sigma_\alpha(M)=\biggl\{ \{q\alpha \} :\,\,\, q \in \Sigma (M)\biggr\}}
$$
will be more convenient. So, as 
$
{\biggl\{ \{a^ub^v\alpha \} :\,\,\, u,v \le M\biggr\}}\subset \Sigma_{\alpha} \left( (ab)^{M}\right)
$,
 Theorem A states that under certain conditions the set $\Sigma_\alpha (M)$ 
 is $\frac{1}{{(\log \log  \log M)}^{\kappa(a,b)}}$-dense in $[0,1]$.

\vskip+0.3cm
For rational $\alpha = \frac{A}{Q}, (A,Q) = 1$ the following result was proven in \cite{BLMV}. 

\vskip+0.3cm

%\textcolor{blue}{Поправил скобки в формулировке. Ещё - может нам переписать формулировку Теоремы A так, чтобы там тоже тройной логарифм возник?}\\

{\bf Theorem B.}\, {\it
Let $a, b$ be multiplicatively independent and  $(ab, Q) = 1$. Then for any $A$ coprime to $Q$ 
the set 
$$
 {
\biggl\{ \{a^ub^v\alpha\}: \,\, u,v \le 3\log Q\biggr\}\,\,\,\,\,
\text{where}\,\,\,\,\, \alpha = \frac{A}{Q}
}
$$
  is  $\frac{\kappa_1(a,b)}{(\log\log\log{Q})^{\kappa_2(a,b)}}$-dense in $[0,1]$ with effective  positive $\kappa_j (a,b), j =1,2$.
}
\vskip+0.3cm

The methods of the proof rely on entropy theory, in particular on an effective version of Rudolph-Johnson's theorem (see \cite{J,R}).

\vskip+0.3cm

In the present paper we would like to give a very  simple and explicit exposition of the effective results from the wonderful paper  \cite{BLMV}, which does not use neither entropy nor measures. Our simplification follows the ideas from \cite{BLMV} but it uses only pigeon hole principle and bounds for simple exponential sums modulo $a^n$.
 For the simplicity of exposition in this paper we consider the case $(a,b) = 1$ only.
We will prove the following theorem.
\vskip+0.3cm

{\bf Theorem 1.}
\, {\it Let $a, b \ge 2, (a,b) = 1$  be positive integers. Then for any $\delta,\varepsilon >0$ for all integers $Q\ge Q_0 (a,b,\delta,\varepsilon)$  and for any 
\begin{equation}\label{alpha}
\alpha = \frac{A}{Q},\,\,\, (A,Q) = 1.
\end{equation}
the
set
$
{
 \Sigma_\alpha (Q^{1+\delta})}
$
is
$
\frac{1}{(\log\log\log Q)^{\frac{1}{8}-\varepsilon}}$-dense in $[0,1]$}.
    
 \vskip+0.3cm
   
{\bf Remark.}\, In Theorem 1 we do not need additional assumption $(ab,Q) = 1$.    In particular, if $Q=(ab)^{r_1}, A=1,$ the statement about density remains true for the set
  $$
  \left\{
  \{a^{u} b^{v}\}:\, -r_1 \le u,v \le  r_2
  \right\}
  $$  
  for large $r_1$ with $ r_2 \ge \delta r_1$. 
  {Also, for $A=1$ Theorem 1 becomes just a statement about density of the set $\Sigma_{\frac{1}{Q}} (Q^{1+\delta})$ and 
  application of the argument from Section 3 below shows that   simply the set $\Sigma_{\frac{1}{Q}} (Q)$ is
  $\frac{1}{(\log Q)^{\kappa_3(a,b)}}$-dense in $[0,1]$  (see (\ref{xi}) with $\kappa_3(a,b) =\frac{1}{\beta-1}-\varepsilon$
  where $\beta=\beta (a,b)>1$ comes from (\ref{feel})).
  
   \vskip+0.3cm

  Dirichlet Theorem states that for any $\alpha\in \mathbb{R}$  and for any positive integer $N$ there exist  coprime integers $A, Q$ satisfying
\begin{equation}\label{diri}
\left| \alpha -\frac{A}{Q} \right| \le \frac{1}{QN},\,\,\,\, 1\le Q\le N.
\end{equation}
If  $\alpha$  is irrational,  then there exist infinitely many rational fractions $\frac{A}{Q}$ satisfying
\begin{equation}\label{che}
\left|\alpha-\frac{A}{Q}\right|\le \frac{1}{Q^2},\,\,\,\, (A,Q) = 1.
\end{equation}
On the other hand, similarly to the result for inhomogeneous approximation,  a famous Chebyshev-Hurwitz-Khintchine's theorem 
(see \cite{Kh}) states that 
for any $\varepsilon>0$, any  irrational $\alpha$  and any real $\beta$ there exist infinitely many positive integers  $q$ such that
$$
||q\alpha - \beta|| < \frac{1-\varepsilon}{q\,\sqrt{5}}.
$$

   \vskip+0.3cm

From Theorem 1 we can easily deduce a Chebyshev-type result concerning inhomogeneous approximation of the form
$$
||q\alpha - \beta||,\,\,\,\,\, q\in \Sigma
$$
(here  $||\cdot || = \min_{x\in \mathbb{Z}} |\cdot - x|$  stands for the distance to the nearest integer). 
If $ q\in \Sigma (Q^{1+\delta})$ then for the fraction $\frac{A}{Q}$ satisfying (\ref{che}) we have
$$
\left|q\alpha  -\frac{Aq}{Q}\right|\le \frac{1}{Q^{1-\delta}}.
$$
Theorem 1 claims that the set  of fractional parts
$$
 \Sigma_\alpha (Q^{1+\delta})=\left\{
 \left\{\frac{Aq}{Q}\right\}:\,\,\, q\in \Sigma, \,\, q\le Q^{1+\delta}\right\}
 $$
 is $\frac{1}{(\log\log\log Q)^{\frac{1}{8}-\varepsilon}}$-dense in $[0,1]$.
So as a corollary of Theorem 1 we immediately deduce the following Chebyshev-type statement.

\vskip+0.3cm

{\bf Theorem 2.}\, {\it Let $a, b \ge 2, (a,b) = 1$  be positive integers. 
Let $\alpha $ be real irrational number and $\beta $ be a real number.
Then  for any $\varepsilon>0$ there exist infinitely many integers $q \in \Sigma$ such that
$$
||q\alpha - \beta ||\le \frac{1}{(\log\log \log q)^{\frac{1}{8}-\varepsilon}}.
$$}

\vskip+0.3cm
Theorem 2 gives  us for any irrational $\alpha$  just an infinite sequence of $q\in \Sigma$ satisfying  the approximation property. Applying Dirichlet's Theorem (\ref{diri}) we can immediately deduce  from Theorem 1 a "uniform"  result analogous to Theorem A, which
deals with numbers $\alpha$ which are not very well approximable by rationals. That is, given any large $N$ we establish the existence of $q\in \Sigma(N)$
with small $ ||q\alpha - \beta||$.

Let $ \psi (t)$ be decreasing to zero as $ t\to\infty$ function such that $ \psi (t) = O\left(\frac{1}{t}\right)$,
and $\Psi (t)$ be the function inverse to the function $ t\mapsto \frac{1}{\psi(t)}$. It is clear that $\Psi (t) = O(t), t\to \infty$.
Suppose that $\alpha \in \mathbb{R}\setminus\mathbb{Q}$  is $\psi$-badly approximable, that is
\begin{equation}\label{bad}
||Q\alpha || \ge \psi (Q),\,\,\,\,\,\, \forall Q\in \mathbb{Z}_+.
\end{equation}
Then obviously Dirichlet Theorem for such an $\alpha$ can be modified as follows.
If irrational $\alpha$ satisfies (\ref{bad}), then 
for any positive integer $N$ there exist  coprime integers $A, Q$ satisfying
\begin{equation}\label{diri1}
\left| \alpha -\frac{A}{Q} \right| \le \frac{1}{QN},\,\,\,\, \Psi (N) \le Q\le N.
\end{equation}
In particular, if $ \psi (t) = t^{-k}, k \ge 1$ then $ \Psi (t) = t^{\frac{1}{k}}$ and
$
\log\log \log \Psi (N) = \log (\log \log N - \log k)
.
$

 \vskip+0.3cm
 
 {\bf Theorem 3.}\, {\it Let $a, b \ge 2, (a,b) = 1$  be positive integers. 
Let $\alpha \in \mathbb{R} $ satisfy (\ref{bad}).
Then for any $\delta,\varepsilon >0$ for all integers $N$ large enough  for  $Q$ defined in (\ref{diri1}) the set
$
{
 \Sigma_\alpha (Q^{1+\delta})}
$
is
$
\frac{1}{(\log\log\log Q)^{\frac{1}{8}-\varepsilon}}$-dense in $[0,1]$, and hence the set
$
 \Sigma_\alpha (N^{1+\delta})
\supset
 \Sigma_\alpha (Q^{1+\delta})
$
is
$
\frac{1}{(\log\log\log \Psi(N))^{\frac{1}{8}-\varepsilon}}$-dense in $[0,1]$.

}.

\vskip+0.3cm

The proof of Theorem 1  will be given in Sections  3--6. In Section 3  we collect  classical results devoted to integer points, distribution of the sequence $s_j$ and application of linear forms in logarithms of integers. Section 4 deals with some combinatorial statements about 
the digital structure of the sets of integer parts $[a^n\{q\alpha\}],\, \, q\in \Sigma (N)$ with respect to the base $a$. Section 5 involves analytic consideration of exponential sums over 
$\mathbb{Z}/{a^\ell}\mathbb{Z}$. In Section 6 we finalise the proof.

\vskip+0.3cm

{\bf 2. Few words about lower bounds.}

\vskip+0.3cm

We would like to cite certain results concerning lower bounds in the problem. 
Let $\Sigma \subset \mathbb{Z}_+$ be written as a sequence
$$
\Sigma : \,\,\, q_1<q_2< q_3<...<q_\nu< q_{\nu+1}<...
$$
of distinct integers in increasing order

In \cite{M} the second author by means of Peres-Schlag's method (see \cite{PS}) proved that the set
$$
\mathcal{B} =
\{ \alpha \in \mathbb{R}: \, 
\inf_{\nu \in
\mathbb{Z}_+ 
}\, \sqrt{\nu} \log \nu \cdot ||q_\nu \alpha|| >0\} =
\{ \alpha \in \mathbb{R}: \, 
\inf_{q \in
\Sigma 
}\, \log q  \log\log  q \cdot  ||q\alpha|| >0\} 
$$
(see formula (\ref{qnu}) from Section 3 below)
has full Hausdorff dimension. Moreover, concerning the problem of distribution of fractional parts of the form
$\{q\alpha\},\,\, q\in \Sigma$ in $[0,1]$, application of Peres-Schlag's method shows that for any sequence 
$$
\pmb{\gamma}:\,\,\, \gamma_\nu\in [0,1],\,\,\,\, \nu =1,2,3,...
$$
 the corresponding set
$$
\mathcal{B}_{\pmb{\gamma}} =
\{ \alpha \in \mathbb{R}: \, 
\inf_{\nu \in
\mathbb{Z}_+ 
}\, \sqrt{\nu} \log \nu \cdot ||q_\nu\alpha -\gamma_\nu|| >0\} 
$$
has also full Hausdorff dimension.
In recent paper \cite{BH} Badziahin and Harrap proved
that  for any $\varepsilon >0$ the set 
$$
 \{ \alpha \in \mathbb{R}: \, 
\inf_{q \in
\Sigma 
}\, (\log q)^{1+\varepsilon} \cdot  ||q\alpha|| >0\} 
$$
is a  Cantor-winning set and has full Hausdorff dimension.
Moreover, they showed that application of the method of the paper \cite{SV} also  proves that the set 
$\mathcal{B}$ has full Hausdorff dimension. As far as we know it is not known if the set $\mathcal{B}$ is a
Cantor-winning set or is winning in some other game.

\vskip+0.3cm

 \vskip+0.3cm

{\bf 3. Set $\Sigma (M)$ and  its difference set.}

%\textcolor{red}{Please explain this inequality. Do I understand correctly that $\beta$ is an irrationality measure of the number $\log a/\log b$?}

\vskip+0.3cm

For our purpose we need a lower bound of the form
\begin{equation}\label{feel}
\left|\frac{\log a}{\log b} - 
\frac{p}{q} \right|> \frac{c}{q^{\beta}},\,\,\,\,\,\ \text{for all rational}\,\,\,\ \frac{p}{q}
\end{equation}
with effective $ c>0, \beta>2$ depending on $a$ and $b$.
Such a bound was firstly obtained by Feldman
\cite{F1,F2} based on a breakthrough result by Baker \cite{Ba}.
Here we would like to refer to recent book \cite{BU} which contains a lot of information about the history and applications of linear forms of logarithms of algebraic numbers. For $ a=2, b=3$ the best known value  $\beta = 5.116201$ was obtained recently by Bondareva, Luchin and Salikhov   by methods of Transcendence Theory (see more general Theorem 1 from \cite{BLS}).

 \vskip+0.3cm
 The number $T(t)$ of integer points in the right angled triangle
$$
x,y >0,\,\,\,\,\,\ x\log a +y \log b \le t
$$
has the asymptotics
\begin{equation}\label{HL}
T(t) = \frac{t^2}{2\, \log a\, \log b} - t\cdot \left(\frac{1}{2\log a}+\frac{1}{2\log b} \right)  + O_{c,\beta}\left(t^{1-\frac{1}{\beta-1}}\right),\,\,\, t \to \infty
\end{equation}
(see 
\cite{HL}, {Theorem 5}).
From (\ref{HL}) we see that for any $t$  large enough and $r_t $  satisfying
$r_t \cdot t^{\frac{1}{\beta-1}} \to \infty$ 
there exist positive integers $u,v$ such that 
\begin{equation}
\label{logabsum}
 t< u\log a + v \log b \le t+r_t.
 \end{equation}
 {Indeed, it is enough to show that $T(t+r_t)-T(t)\to\infty$ as $t\to\infty$.} Or, in other words, {taking exponent of (\ref{logabsum}) we see that} for any $ \tau =e^t$ there exists $q = a^ub^v\in \Sigma \left( \frac{t}{\log a}\right)$ with
 $$
 \tau \le q \le \tau e^{r_t} = \tau (1+ O(r_t)) = \tau + O\left( \frac{\tau}{(\log \tau)^{\frac{1}{\beta-1}-\varepsilon}}\right)
 $$
 for any positive $\varepsilon$. {Note that  from (\ref{HL})  the asymptotic equalities
 \begin{equation}\label{HLA}
 |\Sigma(M) | = T(\log M) \sim \frac{\log^2M}{2\log a \log b},\,\,\,\,\, M \to \infty
 \end{equation}
 and 
 \begin{equation}\label{qnu}
 \log q_\nu \sim \sqrt{2\nu\log a \log b},\,\,\,\,\, \nu \to \infty
 \end{equation}
 follow,
 as well as the inequality
 \begin{equation}\label{HL1}
 q_{\nu+1}- q_\nu \ll_{a,b} \frac{q_\nu}{(\log{q_\nu})^{\frac{1}{\beta-1}-\varepsilon}}.
 \end{equation}
 }

%\textcolor{yellow}{дальше тут была ПЛОХАЯ ошибка, максимум был взят не там. Пришлось исправлять. Я сейчас переделал:}
%\textcolor{blue}{Denote 
%$$
%D_d=\max\limits_{0\le i\le k} \bigl(q_{\nu+i+1}-q_{\nu+i}\bigr).
%$$
%}
%Then \textcolor{blue}{for any $j\le k$ one has}
%\begin{equation}\label{xi}
%\xi_{j+1} - \xi_j \le D_d \asymp  \frac{d}{(\log{d})^{\frac{1}{\beta-1}-\varepsilon}}.
%\end{equation}\vskip+0.3cm

Now we formulate a corollary from the classical results mentioned above.
 \vskip+0.3cm
%\textcolor{red}{(Lemma and the argunemt after it is rewritten using the notation $\Xi(M)$)}

{
{\bf Lemma 1.}\, {\it    Let $\alpha$ be defined in (\ref{alpha}). 
For some positive integers
$M$ and $M_1$, satisfying
\begin{equation}\label{em}
 M <Q,\,\,\,\,\, M_1 = MQ
 \end{equation}
 the set 
$$
\Bigl(\Sigma_\alpha(M_1)-\Sigma_\alpha(M_1)\Bigr)\cup\Bigl(1-\bigl(\Sigma_\alpha(M_1)-\Sigma_\alpha(M_1)\bigr)\Bigr)
  $$
   is $\Delta$-dense  in $[0,1]$ with
   \begin{equation}\label{ed}
   \Delta \ll_{a,b} 
   \frac{1}{(\log\log{M})^{\frac{1}{\beta-1}-\varepsilon}}.
\end{equation}
   .}}

 \vskip+0.3cm
 
 {
 Proof.  
 From (\ref{em}) we see that all the numbers $q \in \Sigma (M)$ are different modulo $ Q$.  
 So all the fractional parts $\left\{ \frac{Aq}{Q}\right\}, \, q\in \Sigma (Q)$ are also different.
 From (\ref{HLA})  by pigeonhole principle we see that there exist different
 $\eta ' =\{ q'\alpha\},\eta''=\{q''\alpha\} \in \Sigma_\alpha (M)$ such that 
 $$
 \frac{1}{Q} \le \eta'-\eta''\le \frac{1}{|\Sigma (M)|} \asymp_{a,b} \frac{1}{\log^2 M}.
  $$
  We define $ d$ by 
  \begin{equation}\label{de}
  \frac{1}{d} = \eta'-\eta'',\,\,\,\, |\Sigma (M)| \le d\le Q.
  \end{equation}
  and consider the collection 
  of all the elements
$q\in \Sigma (d) $  enumerated in the increasing order
\begin{equation}\label{collection}
 q_{1}<q_{2}<...<q_{i}<q_{i+1}<...<q_{k}\le d <q_{k+1},\,\,\,\, q_j \in \Sigma.
\end{equation}
%\textcolor{red}{Я не понимаю, а нужна ли нам оценка снизу в $\sqrt{d}$? Мы же фактически смотрим только на оценку сверху величины $q_{i+1}-q_i$, может сразу с $1$ $q_i$ брать?}

  Denote
$$
D_d=\max\limits_{1\le i\le k} \bigl(q_{i+1}-q_{i}\bigr).
$$
Then  by (\ref{HL1}) for any $i\le k$ one has
\begin{equation}\label{xi}
q_{i+1} - q_{i} \le D_d \ll_{a,b}  \frac{d}{(\log{d})^{\frac{1}{\beta-1}-\varepsilon}}.
\end{equation}
}

 {
As $ d\le Q$,  for $q_j$ from the sequence (\ref{collection}) we have  $q_j q', q_j q'' \in \Sigma(MQ)= \Sigma(M_1)$. 
{Denote}
  $$
  {\eta_j = \frac{q_j}{d}= q_j(\eta'-\eta'')=q_j(\{q'\alpha\}-\{q''\alpha\})=
  \begin{cases} 
  \{q_jq'\alpha\}-\{q_jq''\alpha\} \in \Sigma_\alpha\left( M_1\right)-\Sigma_\alpha \left( M_1\right) & \text{if}\ \ \{q_jq'\alpha\}>\{q_jq''\alpha\} \\
  1-(\{q_jq'\alpha\}-\{q_jq''\alpha\}) \in 1-(\Sigma_\alpha\left( M_1\right)-\Sigma_\alpha \left( M_1\right)) & \text{otherwise}.
  \end{cases}
}
  $$
  Taking into account the estimate (\ref{xi}), we obtain
  $$
 {\eta_{j+1}-\eta_j}=(q_{j+1}-q_j)(\eta'-\eta'') = \frac{q_{j+1}-q_j}{d}
  \le \frac{D_d}{d}
  .
  $$  
  Now we put
  $$
  \Delta = \frac{D_d}{d}\ll_{a,b} \frac{1}{(\log{d})^{\frac{1}{\beta-1}-\varepsilon}} \le \frac{1}{(\log{|\Sigma(M)|})^{\frac{1}{\beta-1}-\varepsilon}}
  \ll_{a,b} \frac{1}{(\log\log{M})^{\frac{1}{\beta-1}-\varepsilon}}.
  $$
 {One can easily see that the union of the sequences $\eta_j$ and $1-\eta_j$, $1\le j \le  k$} forms a $\Delta$-net in the set $ [0,1]$.
   Lemma is proven.$\Box$
   }

 \vskip+0.3cm
 Now we take {
 \begin{equation}\label{eni}
 N = a^{n},\,\
 \,\,\,\,\,
 N = C (\log \log{M})^{\frac{1}{\beta-1}-\varepsilon},\,\,\,
 n \sim  \frac{\frac{1}{\beta -1}-\varepsilon}{\log a} \cdot \log\log\log M
 \end{equation}
 with $C=C(a,b)$  such that $N$ satisfies the inequality
  \begin{equation}\label{eni1}
N\le \frac{1}{\Delta}.
 \end{equation}
 By (\ref{ed}), one can see that the parameters $n$ and $N$, satisfying (\ref{eni}) and (\ref{eni1}), exist.
 }
 Then we consider integer parts
 $$
 [a^n\eta],\,\,\,\,\, \eta \in  \Sigma_\alpha (M_1).
 $$
 By (\ref{eni1})  and the conclusion of Lemma 1 we see that 
 the set
   $$
  (a^{n}\eta'-a^{n}\eta'')\cup(a^n-(a^{n}\eta'-a^{n}\eta'')),\,\,\,\,\,\eta',\eta'' \in \Sigma_\alpha(M_1)
   $$
   is 1-dense in $[0,N]$. 
   So by the inequality
   $$
   |[x-y]-([x]-[y])|\le 1
   $$
   the set of differences 
   $$
     [a^{n}\eta']-[a^{n}\eta''],\,\,\,\,\,\eta',\eta'' \in \Sigma_\alpha (M_1)
   $$
   has at least $N/4$ different elements in $[0,N]$. 
This means that the cardinality of the set
$$
\{
   [a^{n}\eta],\,\,\,\,\,\eta \in \Sigma_\alpha (M_1)\}\cap [0,N]
   $$
   must be greater than $\frac{\sqrt{N}}{2}$. We formulate this conclusion as following\\

  \vskip+0.3cm {
  {\bf Lemma 2.}\,{\it
  Given a positive integer $R$, consider the set
 \begin{equation}\label{theset}
 \frak{X}_n^{R} =\left\{ x\in \mathbb{Z}_+: \, 0\le x \le a^n-1,\,\,\exists \eta \in  {\Sigma_\alpha (R)}\,\,
\text{ such that}\,\, \frac{x}{a^n}\le \eta<\frac{x+1}{a^n}\right\} =
\end{equation}
$$
=
\{\text{\rm numbers  from  the interval}\,\,\, [0,N]\,\,\,\text{\rm of the form}\,\,[a^n\eta],\, \, \eta \in \Sigma_\alpha (R)\}.
$$
Then for $R=M_1$ for the  cardinality $X_n^{M_1}$ of the set   $\frak{X}_n^{M_1} $ we have a lower bound
\begin{equation}\label{comki}
X_n ^{M_1}\ge 
  \frac{\sqrt{N}}{2}.
  \end{equation}
  }
  }

  \vskip+0.3cm

{\bf 4. Combinatorial part.}

\vskip+0.3cm

Let
  $ s\le n$. Together with the set $\frak{X}_n^R$  defined in (\ref{theset}) we consider the set
$$
\frak{X}_{n,s}^R = \{ x\in \mathbb{Z}_+:\,\, 0\le x \le a^s-1,\,\,\, \exists x^*\in \frak{X}_n^R\,\,\text{such that}\,\, x\equiv x^*\pmod{a^s}\}.
$$

%\textcolor{red}{I disagree with the inclusion $\frak{X}_{n,s}\subset  \frak{X}_{s}$. As long as you omit the depndence on $M_1$, it at least should be the same as for the subset as for the bigger set. From the proof I can see that $\frak{X}^{M_1}_{n,s}\subset $\frak{X}^{M_1+n-s}_{s}$. So I suggest to alternate the formulation.}

%\textcolor{red}{Фраза ''Moreover if for $x^*\in \frak{X}_{n}^{M_1}$ there exists $ \eta \in \Sigma_\alpha (M_1)$ satisfying $ \frac{x^*}{a^n}\le \eta <\frac{x^*+1}{a^n}$'' означает, что это некоторое дополнительное предположение. Но разве из того, что $x^*\in \frak{X}_{n}^{M_1}$ не следует автоматически того факта, что найдётся соответствующее $\eta$? В общем, надо как-то формулировку изменить, вам не кажется?}
\vskip+0.3cm

{
{\bf Lemma 3.} \, {\it  Let 
\begin{equation}\label{m2}
M_2 = M_1\cdot  a^{n-s}.
\end{equation}
 If for $x^*\in \frak{X}_{n}^{M_1}$ there exists
$ \eta \in \Sigma_\alpha (M_1)$ satisfying $ \frac{x^*}{a^n}\le \eta <\frac{x^*+1}{a^n}$,
then for $ x\in \frak{X}_{n,s}^{M_2}$ with $ x\equiv x^*\pmod{a^s}$ there exists
$ \eta_1 \in \Sigma_\alpha (M_2)$ satisfying $ \frac{x}{a^s}\le \eta_1 <\frac{x+1}{a^s}$.
In particular, 
$\frak{X}_{n,s}^{M_1}\subset  \frak{X}_{s}^{M_2}$.
}
}

\vskip+0.3cm

Proof.
Let 
$$
\frac{x^*}{a^n} \le \eta < \frac{x^*+1}{a^n},\,\,\,\, \eta = \{a^ub^v\alpha\}  \in \Sigma_\alpha(M_1).
$$
We multiply these inequalities by $a^{n-s}$. Then
$$
\frac{x^*}{a^s} \le a^{n-s}\eta < \frac{x^*+1}{a^s}.
$$
We have
$x^* = x +\lambda a^s, \lambda\in \mathbb{Z}$ and so
$$
\frac{x}{a^s}+\lambda \le a^{n-s}\eta < \frac{x}{a^s} +\lambda +\frac{1}{a^s}.
$$
We see that $ [a^{n-s}\eta] = \lambda$ and
{
$$
a^{n-s} \eta-\lambda = \{a^{n-s}\eta\} =\{ a^{n-s}\{a^ub^v\alpha\}\} =\{a^{u+n-s}b^v\alpha\}=\eta_1\in \Sigma_\alpha(M_1\cdot a^{n-s}).
$$
}
So 
$$
\frac{x}{a^s} \le \eta_1 < \frac{x+1}{a^s},
$$
and $ x \in \frak{X}_s^{M_2}$.$\Box$

\vskip+0.3cm

We consider nonnegative integers $ \ell \le s$ and $\lambda$ from the interval $ 0\le \lambda< a^{s-\ell}$. We will deal with the sets
$$
\frak{X}_{n,s,\ell}(\lambda)  = \{  x\in \frak{X}_{n,s}^{{M_1}}\,\,\text{such that}\,\, x\equiv \lambda\pmod{a^{s-\ell}}\}.
$$
So $\frak{X}_{n,s,\ell}(\lambda)$
consists of certain integers $x$ of the form %\textcolor{red}{(as far as I understand,  the first term should be $ \lambda$ instead of $ \lambda\cdot a^{s-\ell -1}$)}
$$
x = {\lambda}+ x_{s-\ell}a^{s-\ell}+...+ x_{s-2} a^{s-2} + x_{s-1}a^{s-1},\,\,\,\,\,\,\,\,
0\le x_j\le a-1, \, s-l \le j \le s-1.
$$
It is clear that 
\begin{equation}\label{uni}
\frak{X}_{n,s}^{{M_1}} = \bigcup_{\lambda=0}^{a^{s-\ell}-1} \frak{X}_{n, s,\ell}(\lambda).
\end{equation}

\vskip+0.3cm
To avoid too cumbersome notation we denote the cardinalities of the sets $ \frak{X}_n^{M_1}, \frak{X}_{n,s}^{M_1},\frak{X}_{n,s,\ell} (\lambda)$ by 
 $X_n, X_{n,s}, X_{n, s,\ell}(\lambda) $ correspondingly. Also in the next lemma we omit the upper index $M_1$ and write simply
 $\frak{X}_n =  \frak{X}_n^{M_1}, \frak{X}_{n,s} =\frak{X}_{n,s}^{M_1}$.

 \vskip+0.3cm
{\bf Lemma 4 (Main Combinatorial Lemma).} 
{\it
Let $\ell = o(n)$. Suppose that $X_n\ge c a^{\frac{n}{2}}$ with some positive $c>0$.
Then for every  $\varepsilon >0$ and $ n\ge n_0 (c, \varepsilon)$ there exists
$s$  from the interval $ \varepsilon n \le s \le n$  and
$\lambda$ from the interval  $0\le \lambda {<a^{s-\ell}}$ such that  for the cardinality of the set $\frak{X}_{n,s,\ell} (\lambda)$  we have lower bound
\begin{equation}\label{slsl}
X_{n,s,\ell} (\lambda) \ge a^{\left(\frac{1}{2}- 2\varepsilon\right)\ell}.
\end{equation}

}

\vskip+0.3cm
Proof. 
We consider integers
$$
n_j = n-j\ell,\,\,\,\, 0\le j \le J,\,\,\,\, J = \left\lceil(1-\varepsilon) \frac{n}{\ell}\right\rceil.
$$
$$
\varepsilon n -\ell < n_J = n - J \ell \le \varepsilon n
$$
and suppose $n$ to be large enough.
As there are just $a^{n_J}$ residues $\pmod{a^{n_J}}$, we can find
 $\lambda_0\pmod{a^{n_J}}$
 such that the set
 \begin{equation}\label{l1}
 \frak{X}'_n = \frak{X}_n\cap \{ x\pmod{a^n}: \,\, x\equiv \lambda_0 \pmod{a^{n_J}}\} \subset \frak{X}_n
 \end{equation}
 has the cardinality
 $$
  X'_n \ge \frac{X_n}{a^{n_J}}
\ge N^{\eta},\,\,\,\, \eta = \frac{1}{2}- 2\varepsilon.
$$
Then we consider the sets
$$ 
\frak{X}'_{n,s}
=
 \{ x'\pmod{a^s}:\,\, \exists x^*\in \frak{X}_n'\,\,\text{such that}\,\, x'\equiv x^*\pmod{a^s}\}
,\,\,\,\,\, s\le n.
$$
Here we should note that by the construction the set $\frak{X}'_{n,n_J}$  consists just of one element and so
\begin{equation}\label{l3}
X'_{n,n_J} =1.
\end{equation}

We shall prove that for some
$j\in 
\{ 0,1,...,J\}
$ 
and
$ \lambda\pmod{a^{n_j-\ell}}$   for the cardinality $X'_{n,n_j,\ell}(\lambda)$ of the set
$$
\frak{X}'_{n,n_j,\ell}(\lambda) =\{ x\in \frak{X}'_{n,n_j}:\,\,\,\,\, x\equiv\lambda \pmod{a^{n_j-\ell}}\}
%\quad\text{\textcolor{red}{Модуль д.б. $a^{n_j-\ell}$, мне кажется}}
$$
the inequality
\begin{equation}\label{l2}
X'_{n,n_j,\ell}(\lambda) \ge  a^{\eta \ell}
\end{equation}
is valid. 
%\textcolor{red}{А разве не должно в таком случае в формулировке леммы быть $a^{\left(\frac{1}{2}-2\varepsilon\right)\ell}$?}
This would be enough to finish the proof of the lemma, because 
$\frak{X}'_{n,n_j,\ell}(\lambda)\subset \frak{X}_{n,n_j,\ell}(\lambda)$, 
%\textcolor{red}{Написана тавтология, видимо в правой части штрих не нужен в обозначении множества.} 
and
we get  (\ref{slsl}) for $ s =n_j$.

So now we are proving (\ref{l2}).
Suppose that  the desired statement is not valid and
\begin{equation}\label{hat}
X'_{n, n_j,\ell}(\lambda) < a^{\eta\ell}
\end{equation}
for every 
$ j \in 
\{ 0,1,..., J\} $
and for every 
$ \lambda_j \pmod{a^{n_j-\ell}}$.
From the equality
\begin{equation}\label{hat1}
\frak{X}_{n,n_j}' = \bigcup_{\lambda \pmod{a^{n_j-\ell}}} \frak{X}_{n,n_j,\ell}'(\lambda).
\end{equation}
(which is analogous to (\ref{uni})) we get
\begin{equation}\label{l4}
X_{n,n_j}' > a^{\eta n_j},\,\,\,\,\, 0\le j \le J.
\end{equation}

Indeed, inequality (\ref{l4}) can be proven by induction.
The base of the induction is verified for 
$j=0, \,\, n_0 = n$ by the inequality
$$
X'_n = X_{n_0}' > N^{\eta} = a^{\eta n_0}.
$$
To prove the inductive step we suppose that 
\begin{equation}\label{hat2}
X_{n,n_j}' > a^{\eta n_j}
\end{equation}
Then for different  $\lambda_1,\lambda_2 \pmod{a^{n_{j+1}}}, n_{j+1}= n_j -\ell$
we have different disjuncted sets $\frak{X}_{n,n_{j},\ell}' (\lambda_j)$. 
For the cardinality of the set 
 $\frak{X}_{n,n_{j+1}}'$  we have the equality
$$
X_{n,n_{j+1}}'  =|\{ \lambda\pmod{a^{n_{j+1}}}:\,\,\, \frak{X}_{n,n_{j},\ell}' (\lambda)\neq \varnothing\}|.
$$
So from 
(\ref{hat1}) and (\ref{hat},\ref{hat2})
we see that 
$$
X_{n,n_{j+1}}'  \ge \frac{X_{n,n_{j}}'}{\max_\lambda X'_{n,n_j,\ell}(\lambda) }\ge a^{\eta n_{j+1}},
$$
and (\ref{l4}) is proven.

Now from (\ref{l4}) for a particular choice  $ j= J$ we get $ X_{n,n_J}' \ge a^{\eta n_J}>1$.
This contradicts to (\ref{l3}). We proved inequality (\ref{l2}).$\Box$

\vskip+0.3cm

We finish this combinatorial section with defining a set $\frak{Y}$  modulo $a^\ell$.
For integer
$$ 
x=x_0+ ax_1, \,\,\, 0\le x_0 \le a-1
$$
we consider $a$-ary shift
$$
 T_a(x) =\left[\frac{x}{a}\right]= \frac{x-x_0}{a} = x_1 \,\,\,
$$
and the set 
$$
\frak{Y} = T_a^{s-\ell} (\frak{X}_{n,s,\ell}(\lambda)) \subset \{ 0,1,2,...,a^\ell-1\}.
$$
The map
$T_a^{s-\ell} $ forms   a bijection $y = T_a^{s-\ell} (x)$  between $\frak{X}_{n,s,\ell}(\lambda)$ and $\frak{Y}$ where
$$
x = \lambda + x_{s-\ell}a^{s-\ell}+
x_{s-\ell+1} a^{s-\ell+1}+
...+ x_{s-2} a^{s-2} + x_{s-1}a^{s-1}
\,\,\,
\longleftrightarrow 
\,\,\,
y= x_{s-\ell}+ x_{s-\ell+1}a+...+ x_{s -2} a^{\ell -2} + x_{s -1}a^{\ell -1}
$$
or
\begin{equation}\label{corr}
\frac{x}{a^s} =  \frac{\lambda}{a^s} + \frac{y}{a^\ell},\,\,\,\,\, x\in \frak{X}_{n,s,\ell}(\lambda),\,\,\, y\in \frak{Y}.
\end{equation}
So for the cardinality $Y$ of the set $\frak{Y}$ we have the equality
\begin{equation}\label{fraa}
Y = |\frak{Y}| =|\frak{X}_{n,s,\ell}(\lambda)| = {X}_{n,s,\ell}(\lambda).
\end{equation}
{
Moreover, if we put $\gamma = \frac{\lambda}{a^s} $
from all the definitions  and Lemma 3 we see that
\begin{equation}\label{defo}
\forall \, y \in \frak{Y}\,\,\,\, \exists q \in \Sigma (M_2)\,\,\,\,\text{such that}\,\,\,\,
\left| \frac{y}{a^\ell} +\gamma - \{q\alpha\}\right| \le \frac{1}{a^s},
\end{equation}
where $M_2$ depends on $M_1= MQ ,n,s$ and is defined in (\ref{m2}).
}

 \vskip+0.3cm

{\bf 5. Analytic part}

\vskip+0.3cm

We consider coprime  positive integers $a,b$ and positive integer $\ell$.
By $S$ we denote the multiplicative order of $b$ modulo $a^\ell$.
So the successive powers $b^0,b^1,..., b^{S-1}$ form a  cyclic subgroup $\frak{S}$ of the multiplicative group modulo $a^\ell$ and
\begin{equation}\label{oord}
S\, | \, \varphi (a^\ell),\,\,\,\,\,
S \le \varphi (a^\ell) < a^\ell.
\end{equation}
It is well known  (e.g. see Lemma 4 from \cite{Sch} and its corollaries)
that there exists $\kappa_1 = \kappa_1(a,b)>0$ such that 
 $S \ge \kappa_1 a^\ell$. Moreover,  $\frak{S}$ 
contains a cyclic subgroup $\frak{S}'$ of the form
$$
\frak{S}\supset \frak{S}' =\{ x\pmod{a^\ell}:\,\,\, x\equiv 1\pmod{a^{\ell_1}}\}
$$
with $ \ell_1 = 
\ell -  \varkappa$, where  $\varkappa =\lceil a^3\log_a b\rceil$ (see the end of the proof of Lemma 3.3 from \cite{BLMV}).

 The following Lemma 5 we take from \cite{BLMV} without changes (see formula (3.5a) from \cite{BLMV}).
It shows that a  lot of exponential sums over $\frak{S}$ vanish.
In fact the argument from \cite{BLMV}  is close here to Lemma 4 and its corollaries from \cite{Sch}.
The next Lemmas 6 and 7 are just statements from Sections 3.6--3.8 from \cite{BLMV} formulated without measures.

\vskip+0.3cm

 {\bf Lemma 5.} \, 
 {\it For all $\ell$ large enough  and for 
 $ \ell_1 = \ell - \varkappa$ we have a bound for the exponential sum
 \begin{equation}\label{a}
\left| \sum_{s\in \frak{S}} e^{2\pi i \frac{ms}{a^\ell}}
 \right|
 =
 \left| \sum_{w=0}^{S-1} e^{2\pi i \frac{mb^w}{a^\ell}}
 \right|
\le 
\begin{cases}
 0,\,\,\,\,\,\,\,\,  m \not\equiv 0 \pmod{a^{\ell_1}},
\cr
S, \,\,\,\,\,\,\,\,  m \equiv 0 \pmod{a^{\ell_1}}.
\end{cases}
\end{equation}
}

 \vskip+0.3cm
 We do not give here a proof of Lemma 5 but supply all other lemmas with proofs.

 \vskip+0.3cm
{\bf Lemma 6.} \,  {\it For  average of the  sums
$$
 \sigma (m) = \sum_{y\in \frak{Y}} e^{2\pi i \ m \left(\frac{y}{a^\ell}+\gamma\right)},
 $$
{where $\gamma\in\mathbb{R}$ is an arbitrary shift,} we have the bound
\begin{equation}\label{a}
\sum_{s\in \frak{S} }|\sigma(ms)|^2 \le w_m SY
\end{equation}
with
$w_m=a^{\varkappa}{{\rm gcd}(a^{\ell_1},m)}$
.}

\vskip+0.3cm

Proof. Let 
$a_m = \frac{a^{\ell_1}}{{\rm gcd}(a^{\ell_1},m)}$ and define
$\chi_{\frak{Y}} (y)$ to be the characteristic function of the set $\frak{Y}\subset \mathbb{Z}/a^\ell \mathbb{Z}$.
 Now by Lemma 5  and the Cauchy-Schwarz inequality  we have
$$
\sum_{s\in \frak{S}} |\sigma (ms)|^2 =
\sum_{y,y'\in \frak{Y}} 
\sum_{s\in \frak{S}} 
e^{2\pi i \frac{ms(y-y')}{a^\ell}} \le
S  \sum_{\begin{array}{c}y,y'\in \frak{Y}: \cr y\equiv y' \pmod{a_m}\end{array}} 1 =
$$
$$
S\sum_{\begin{array}{c}y''\pmod{a^\ell}: \cr y''\equiv 0\pmod{a_m}\end{array}}
\sum_{y\pmod{a^\ell}}
\chi_{\frak{Y}}(y) \chi_{\frak{Y}}(y+y'')
\le
S\left(
\sum_{\begin{array}{c}y''\pmod{a^\ell}: \cr y''\equiv 0\pmod{a_m}\end{array}} 1\right)\cdot
\left(
\sum_{y\pmod{a^\ell}}
(\chi_{\frak{Y}}(y))^2 \right) =  w_m S Y,
$$ 
and (\ref{a}) is proven.$\Box$

\vskip+0.3cm

 Now we consider 1-periodic  smooth function
$f(z): \mathbb{R} \to \mathbb{R}.$
We are interested in the remainder
\begin{equation}\label{remi}
R_{s}[f,\frak{Y}] = 
\frac{1}{Y} \sum_{y\in \frak{Y}} f \left(s \cdot \left(\frac{y}{a^\ell}+\gamma\right)\right) - \int_0^1 f(t)dt
\end{equation}
for approximation to the integral $\int_0^1 f(t)dt
$.

\vskip+0.3cm
{\bf Lemma 7 (Main Analytic Lemma).} \,{\it  For the remainder (\ref{remi}) we have a bound in average
\begin{equation}\label{b}
\frac{1}{S} \sum_{s\in \frak{S}} \left| R_s[f,\frak{Y}]\right|^2
=
\frac{1}{S} \sum_{w=0}^{S-1} \left| R_{b^w}[f,\frak{Y}]\right|^2
 \ll_{a,b} \frac{||f'||_2^2}{Y},
\,\,\,\text{where }\,\,\,
{|| f'||^2_2}=\int_0^1|f'(z)|^2 dz.
\end{equation}
 In particular, there exists positive $s = b^w \in \frak{S}$ such that
$$
\left|
R_{s}[f,\frak{Y}] \right|\ll_{a,b} \frac{|| f'||_2}{\sqrt{Y}}.
$$
}

\vskip+0.3cm

Proof.
Decomposing $f$ into Fourier series we get
$$
f(t) = \sum_{m\in \mathbb{Z}} f_m e^{2\pi i mt}
,\,\,\,\,\,\text{and}\,\,\,\,\,
R_{s}[f,\frak{Y}] = 
 \frac{1}{Y}  \sum_{m\neq 0} f_m \sigma (ms)
 .
 $$
 Then    by Cauchy-Schwarz  inequality,

 $$
 |
 R_{s}[f,\frak{Y}] 
 |^2
 \le
  \frac{1}{Y^2}  
  \sum_{m\neq 0} (m\, |f_m|)^2 \cdot
  \sum_{m\neq 0} \left(\frac{|\sigma(ms)|}{m}\right)^2 =
    \frac{||f'||_2^2}{Y^2}  
 \cdot
  \sum_{m\neq 0} \left(\frac{|\sigma(ms)|}{m}\right)^2
 ,
\,\,\,\,\,
||f'||_2^2 =
  \sum_{m\neq 0} (m\, |f_m|)^2 
  $$
  Finally, again  by Cauchy-Schwarz we obtain
   $$
  \sum_{s\in \frak{S}} |
 R_{s}[f,\frak{Y}] 
 |^2
    \le \frac{||f'||_2^2}{ Y^2}  \cdot
 \sum_{m\neq 0} \frac{
 1}{m^2} \sum_{s\in \frak{S} }|\sigma(ms)|^2 
     \le \frac{S ||f'||_2^2}{ Y}  \cdot
 \sum_{m\neq 0} \frac{
 w_m}{m^2}
  \ll_{a,b}\frac{S||f'||_2^2}{Y}
  $$
  as
$$
\sum_{m\neq 0}\frac{w_m}{m^2}\ll_a \sum_{m=1}^\infty \frac{{\rm gcd} (a^\ell,m)}{m^2} =
\sum_{d|a^\ell}d \left( \sum_{m: \, d|m}\frac{1}{m^2}\right) \ll
    \sum_{d|a^\ell}\frac{1}{d}\le \frac{a}{\varphi (a)}
 .
 $$
  Inequality (\ref{b}) is proven.$\Box$

\vskip+0.3cm

{\bf Lemma 8.}\,{\it
Suppose that 
\begin{equation}\label{HH}
H=o \left(Y^{\frac{1}{4}}\right) = o\left(X_{n,s,\ell}(\lambda)^{\frac{1}{4}}\right).
\end{equation}
 Then for any $ z\in [0,1]$ there exist integer $ w$ from the interval
 $ 0\le w\le S-1$ and $ x \in \frak{X}_{n,s,\ell}(\lambda)$ such that 
 $$
\left|\left| b^w\frac{x}{a^s} -z\right|\right| \le \frac{1}{H},
 $$
 } 
 
 \vskip+0.3cm

Following  Lemma 4 from \cite{Vino}  we  define a standard function 
$$
f(t) =\sum_{m\in \mathbb{Z}} f_m e^{2\pi i m t}
$$
with

({\bf i}) \, $f(t) = 0 $ if $ ||t||\ge \frac{1}{H}$,

({\bf ii}) \, $f_0 = \frac{1}{H} $,

({\bf iii}) \, $f_m \ll \min \left( \frac{1}{H}, \frac{H}{|m|^{2}}\right) $  for  $m\neq 0$. 

\vskip+0.3cm

Then for any $ z\in [0,1)$ the function
$ f_z(t) = f(t-z)$ satisfies
$$
||f'_z||_2^2 =||f'||_2^2 =
  \sum_{m\neq 0} (m\, |f_m|)^2 
  \ll
  \sum_{0<|m|\le H} \frac{|m|^2}{H^2} + 
    \sum_{|m|\ge H} \frac{H^{2}}{|m|^{2}}
    \ll H.
$$ 

Application of the Main Analytic Lemma for  any $z\in [0,1)$ gives $ w$ from the range $ 0\le w \le S-1$
such that 
$$
\sum_{y\in \frak{Y} }f_z \left(b^w\left(\frac{y}{a^\ell} - \gamma\right)\right) = \frac{Y}{H} + 
O(  YR_{s}[f,\frak{Y}] ) =
\frac{Y}{H} + O\left({H}{\sqrt{Y}}\right)>0
$$
by the choice of $H$ in (\ref{HH}). 
So for any $z$ with some $ w$ from the range $ 0\le w \le S-1$
and $y \in \frak{Y}$ one has 
$$
\left|\left| b^w\frac{x}{a^s} -z\right|\right|=
\left|\left| b^w\left(\frac{y}{a^\ell}-\gamma\right) -z\right|\right|<\frac{1}{H},
$$
where $x\in \frak{X}_{n, s,\ell}(\lambda)$ and $y\in \frak{Y}$ are in correspondence (\ref{corr}).
Lemma is proven.$\Box$

\vskip+0.3cm

 {\bf 6. Proof of Theorem 1.}
 
 \vskip+0.3cm

{
 Recall that by our choice of parameters(\ref{em},\ref{ed},\ref{eni}) we have
 $$
 M
 \le Q,\,\,
 M_1= MQ,\,\,M_2 = M_1\cdot a^{n-s} \le M_1 N,\,\,
 \Delta \ll \frac{1}{(\log\log M)^{\frac{1}{\beta-1}-\varepsilon}},\,\,
 N \asymp (\log\log M)^{\frac{1}{\beta-1}-\varepsilon},\,\,
n\asymp \log N \asymp \log\log \log M.
$$
Now we put
$$
M = Q^{\frac{\delta}{2}}
,\,\,\,\,\,
  \ell =\left[\frac{\log (\frac{s\log a}{2\log b })}{\log a}\right].
  $$
  }
  Then
  $a^\ell \asymp  s \asymp n $
  and
  $ b^{a^\ell}\le  a^{\frac{s}{2}}$.
  
  From (\ref{corr}) and (\ref{defo}) we see that 
  
    \begin{equation}\label{w1}
  \left|\frac{x}{a^s} -\{q\alpha\}\right|< \frac{1}{a^s},\,\,\,\,\text{for some}\,\,\,\,\,
  q \in \Sigma(M_2)
  \end{equation}

  From the conclusion (\ref{comki}) of Lemma 2 we see that condition 
  (\ref{slsl}) of Lemma 4 is satisfied. By (\ref{slsl}) and (\ref{fraa})  we deduce 
  \begin{equation}\label{YY}
  Y \ge c a^{\left(\frac{1}{2}-\varepsilon\right)\ell}\asymp n^{\frac{1}{2}-\varepsilon}\asymp
  (\log\log\log M)^{\frac{1}{2}-\varepsilon}
  ,
  \end{equation}
  by the choice of parameters.
  
 For any $ z\in [0,1]$  Lemma 8  gives an integer $ w$ from the interval
 $ 0\le w\le S-1$ and $ x \in \frak{X}_{n,s,\ell}(\lambda)$ such that 
 
    \begin{equation}\label{w2}
\left|\left| b^w\frac{x}{a^s} -z\right|\right| \le \frac{1}{H},\,\,\,\,\text{where}\,\,\,\,\,
H =Y^{\frac{1}{4}-\varepsilon}.
 \end{equation}
We should not that $ b^w \le b^{a^\ell}\le  a^{\frac{s}{2}}$. So from (\ref{w1}) we get
  $$
  \left|b^w\frac{x}{a^s} -b^w\{q\alpha\}\right|< \frac{b^w}{a^s}\le \frac{1}{a^{\frac{s}{2}}}.
$$
Note that  $
b^w\{q\alpha\} \in [0,1)$ and
$
b^w\{q\alpha\}  = \{
b^w\{q \alpha\} \} =  
\{b^wq\alpha\} $. So,
taking into account (\ref{w1},\ref{YY},\ref{w2}) we see that 
for any $ z\in [0,1]$ there exist 
$q_* =b^wq \in \Sigma (M_1\cdot a^{n-s} \cdot  b^{a^\ell}) $
such that
$$
\left|
 \{q_*\alpha\} - z
 \right| \le  \frac{1}{H} + \frac{1}{a^{\frac{s}{2}}} \asymp
  \frac{1}{(\log\log\log M)^{\frac{1}{8}-\varepsilon}}.
 $$
 Now we should note that by the choice of parameters we have
 $$
 M_1\cdot a^{n - s}\cdot b^{ a^\ell}\le
 MQN^2 \le Q^{1+\delta}
 $$
 for $Q$ large enough.
 So we see that for $Q$ large enough   for  $\alpha =\frac{A}{Q}$ the
set
$
\{ \{q\alpha\}:\,\,\, q \in \Sigma (Q^{1+\delta})\}
$
is
$
\frac{1}{(\log\log\log Q)^{\frac{1}{8}-\varepsilon}}$-dense.$\Box$
  
\vskip+0.3cm
\textbf{Acknowledgements}. 
The work of the second named author  has received funding from the European Research Council (ERC) under the European Union’s Horizon 2020 Research and Innovation Program, Grant agreement no. 754475.
Both authors are supported by the Foundation for the Advancement of Theoretical Physics and Mathematics “BASIS”.
\vskip+0.3cm

\noindent Dmitry Gayfulin,\\
Graz University of Technology, Institute of Analysis and Number Theory,\\
Steyrergasse 30/II, 8010 Graz, Austria\\
and\\
Institute for Information Transmission Problems,\\
19 Bolshoy Karetnyi side-str., Moscow 127994, Russia.\\
\textit{gamak.57.msk@gmail.com}
\vskip+0.5cm
\noindent Nikolay Moshchevitin,\\
Israel Institute of Technology (Technion), Center for Mathematical Sciences\\
and\\
Institute for Information Transmission Problems,\\
19 Bolshoy Karetnyi side-str., Moscow 127994, Russia.\\
\textit{moshchevitin@technion.ac.il, nikolaus.moshchevitin@gmail.com}

\end{document}